\documentclass[a4paper,11pt]{article}
\usepackage[T1]{fontenc}
\usepackage[cp1250]{inputenc}
\usepackage{lmodern}

\usepackage[english]{babel}
\usepackage{latexsym}
\usepackage{amsmath}
\usepackage{amsthm}
\usepackage{amssymb}
\usepackage{amsfonts}
\usepackage{stackrel}
\usepackage{dsfont}
\usepackage{tikz}
\usepackage{slashbox}
\usepackage[mathscr]{eucal}
\usepackage{authblk}

\newtheorem{thm}{Theorem}

\newenvironment{pro}{\begin{flushleft} \textbf{Proof}\\* \end{flushleft}}{\hfill\(\blacksquare\) \\ }

\newcommand{\reals}{\mathbb{R}}
\newcommand{\naturals}{\mathbb{N}}

\newcommand{\eps}{\varepsilon}

\newcommand{\Ffamily}{\mathcal{F}}

\newcommand{\Vfamily}{\mathcal{V}}
\newcommand{\Wfamily}{\mathcal{W}}

\newcommand{\vphi}{\varphi}
\newcommand{\DST}{\mathcal{D}}

\newcommand{\Addresses}{{
  \bigskip
  \footnotesize

  Krukowski M. (corresponding author), \textsc{Technical Univeristy of \L\'od\'z, \ Institute of Mathematics, \ W\'ol\-cza\'n\-ska 215, \
90-924 \ \L\'od\'z, \ Poland}\par\nopagebreak
  \textit{E-mail address} : \texttt{krukowski.mateusz13@gmail.com}

  \medskip

}}

\begin{document}
\title{Arzel\`a-Ascoli theorem in uniform spaces}
\author{Mateusz Krukowski}
\affil{Technical University of \L\'od\'z, Institute of Mathematics, \\ W\'ol\-cza\'n\-ska 215, \
90-924 \ \L\'od\'z, \ Poland}
\maketitle

\begin{abstract}
In the paper, we generalize the Arzel\`a-Ascoli theorem in the setting of uniform spaces. At first, we recall well-known facts and theorems coming from monographs of Kelley and Willard. The main part of the paper introduces the notion of extension property which, similarly as equicontinuity, equates different topologies on $C(X,Y)$. This property enables us to prove the Arzel\`a-Ascoli theorem for uniform convergence. The paper culminates with applications, which are motivated by Schwartz's distribution theory. Using the Banach-Alaoglu-Bourbaki theorem, we establish relative compactness of subfamily of $C(\reals,\DST'(\reals^n))$.
\end{abstract}

\smallskip
\noindent 
\textbf{Keywords : } Arzel\`a-Ascoli theorem, uniform space, uniformity of uniform convergence (on compacta)

\section{Introduction}

\noindent
Around 1883, Cesare Arzel\`a and Giulio Ascoli provided a necessary and sufficient conditions under which every sequence of a given family of real-valued continuous functions, defined on a closed and bounded interval, has a uniformly convergent subsequence. Since then, numerous generalizations of this result have been obtained. For instance, in \cite{Munkres} on page 278, the compact subsets of $C(X,\reals)$, with $X$ a compact space, are exactly those which are equibounded and equicontinuous. The space $C(X,\reals)$ is given the standard norm 

$$\| f\| := \sup_{x\in X} |f(x)|$$

\noindent
With a change of norm topology to the topology of uniform convergence on compacta, one can obtain a version of Arzel\`a-Ascoli theorem for locally compact space $X$ and metric space $Y$ (\cite{Munkres} on page 290). The author, together with prof. Bogdan Przeradzki endeavoured to retain the topology of uniform convergence (comp. \cite{KrukowskiPrzeradzki}). To this end, the $BC(X,Y)$-extension property was employed. \\
In this paper, we improve on the result obtained by \cite{KrukowskiPrzeradzki} by considering the uniform spaces, which are thoroughly described in \cite{James} in chapter 7 and 8 (beginning on page 103). The concept generalizes metric spaces and provides a convenient setting for studying uniform continuity and uniform convergence. A brief recap on uniformity of uniform convergence and uniformity of uniform convergence on compacta is presented in section \ref{sectionrecap}. We highlight the importance of equicontinuity, which equates the topology of pointwise convergence, the compact-open topology and the topology of uniform convergence on compacta. This property of 'bringing the topologies together' is a motivation for the extension property, whose objective is to equate the topology of uniform convergence with the aforementioned ones. \\
The main part of the paper is section \ref{sectionmainpart}, where we build on the Arzel\`a-Ascoli theorem in \cite{Kelley} on page 236 (theorem \ref{arzela-ascoliforuniformconvergenceoncompacta} in this paper). We introduce a new concept of extension property and investigate how it relates to $BC(X,Y)$-extension property studied in \cite{KrukowskiPrzeradzki}. The culminating point is theorem \ref{ArzelaAscoli}, which characterizes compact subsets of $C(X,Y)$ with the topology of uniform convergence (rather than uniform convergence on compacta). The space $X$ is assumed to be locally compact and $Y$ is a Hausdorff uniform space. Finally, in section \ref{sectionapplications}, we present two possible applications of theorem \ref{ArzelaAscoli} in the theory of distributions. The space $\DST'(\reals^n)$ is known to be nonmetrizable (\cite{Brezis} on page 81), yet it can be given a uniform structure.

\section{Uniform convergence and uniform convergence on compacta}
\label{sectionrecap}

\noindent
We recall the concept of uniformity of uniform convergence. A comprehensive study of this notion can be found in \cite{Kelley} on page 226 or \cite{Willard} on page 280. For convenience, we extract the essence of these discussions in the following theorem:

\begin{thm}(uniformity of uniform convergence)\\
Let $X$ be a set, $(Y,\Vfamily)$ be a uniform space, $\Ffamily \subset Y^X$ and define $\dagger : \Vfamily \rightarrow 2^{\Ffamily \times \Ffamily}$ by

\begin{gather}
\forall_{V \in \Vfamily} \ V^{\dagger} := \bigg\{ (f,g) \in \Ffamily \times \Ffamily \ : \ \forall_{x \in X} \ (f(x),g(x))\in V \bigg\}
\label{basetuc}
\end{gather}

\noindent
Then the family $\{V^{\dagger} \ : \ V \in \Vfamily\}$ is a base for uniformity $\Wfamily$ for $\Ffamily$, which we call the uniformity of uniform convergence. 
\end{thm}

\noindent
The topology induced by the uniformity of uniform convergence is called the \textit{topology of uniform convergence}, which we denote by $\tau_{uc}$. The base sets for this topology are of the form

\begin{gather}
V^{\dagger}[f] = \bigg\{ g \in \Ffamily \ : \ \forall_{x \in X} \ g(x) \in V[f(x)] \bigg\}
\label{dagger}
\end{gather}

\noindent
where $V[y] = \{z \in Y \ : \ (y,z) \in V\}$ as in \cite{Kelley} on page 176. A concept, which is closely related to the uniformity of uniform convergence is the \textit{uniformity of uniform convergence on compacta}. It appears in \cite{Kelley} on page 229 or \cite{Willard} on page 283.

\begin{thm}(uniformity of uniform convergence on compacta)\\
Let $(X,\tau_X)$ be a topological space, $(Y,\Vfamily)$ be a uniform space and $\Ffamily \subset Y^X$. The family 

\begin{gather}
\bigg\{ (f,g) \in \Ffamily \times \Ffamily \ : \ \forall_{x \in D} \ (f(x),g(x))\in V\bigg\}_{V \in \Vfamily, \ D \Subset X}
\label{subbasetucc}
\end{gather}

\noindent
is a subbase for uniformity $\Wfamily$ for $\Ffamily$, which we call the uniformity of uniform convergence on compacta. \\
\textbf{Remark : } The notation $D \Subset X$ stands for '$D$ is a compact subset of $X$'.
\label{uniformityuniformconvergenceoncompacta}
\end{thm}

\noindent
The topology induced by the uniformity of uniform convergence is called the \textit{topology of uniform convergence on compacta}, which we denote by $\tau_{ucc}$. It is an easy observation that $\tau_{ucc}$ is weaker than $\tau_{uc}$. In fact, if we restrict our attention to $C(X,Y)$ rather than $Y^X$, which we do in the sequel, $\tau_{ucc}$ becomes a familiar compact-open topology (proof in \cite{Kelley} on page 230 or \cite{Willard} on page 284).\\
Another important concept is that of \textit{equicontinuity} (\cite{Kelley} on page 232 or \cite{Willard} on page 286). For a topological space $(X,\tau_X)$ and a uniform space $(Y,\Vfamily)$, we say that the family $\Ffamily \subset Y^X$ is equicontinuous at $x \in X$ if 

\begin{gather}
\forall_{V \in \Vfamily} \ \exists_{U_x \in \tau_X} \ \forall_{f \in \Ffamily} \ f(U_x) \subset V[f(x)]
\label{equi}
\end{gather}

\noindent
The family $\Ffamily$ is said to be equicontinuous if it is equicontinuous at every $x \in X$. If this is the case, then the topology of pointwise convergence coincides with the compact-open topology, i.e. $\tau_{pc}|_{\Ffamily} = \tau_{co}|_{\Ffamily}$ (proof in \cite{Kelley} on page 232 or \cite{Willard} on page 286). Together with what we said above and the fact that equicontinuous family is in particular continuous, we obtain that the topology of pointwise convergence coincides with the topology of uniform convergence on compacta on $\Ffamily$, i.e. $\tau_{pc}|_{\Ffamily} = \tau_{ucc}|_{\Ffamily}$.

\section{Arzel\`a-Ascoli theorems}
\label{sectionmainpart}

\noindent
Our objective in this section is to present the Arzel\`a-Ascoli theorem for the topology of uniform convergence in the setting of uniform spaces. Our starting point is the theorem provided in \cite{Kelley} on page 236.

\begin{thm}(Arzela-Ascoli for uniform convergence on compacta)\\
Let $(X,\tau_X)$ be a locally compact space, $(Y,\Vfamily)$ be a Hausdorff uniform space and let $C(X,Y)$ have the topology of uniform convergence on compacta $\tau_{ucc}$. A subfamily $\Ffamily \subset C(X,Y)$ is relatively $\tau_{ucc}$-compact if and only if 

\begin{description}
	\item[\hspace{0.4cm} (AAucc1)] $\Ffamily$ is pointwise relatively compact, i.e. $\{f(x) \ : \ f \in \Ffamily\}$ is relatively $\tau_Y$-compact at every $x \in X$ 
	\item[\hspace{0.4cm} (AAucc2)] $\Ffamily$ is equicontinuous
\end{description}  
\label{arzela-ascoliforuniformconvergenceoncompacta}
\end{thm}

\noindent
Our line of attack is motivated by equicontinuity, which equates three topologies, of pointwise convergence, compact-open and uniform convergence on compacta. We look for a property that will do the same thing with uniform convergence. \\
For a topological space $(X,\tau_X)$, uniform space $(Y,\Vfamily)$ and $\Ffamily \subset Y^X$ suppose that $\Wfamily_{ucc}$ is a uniformity of uniform convergence on compacta on $\Ffamily$. We say that $\Ffamily$ satisfies the \textit{extension property} if 

\begin{gather}
\forall_{V \in \Vfamily} \ \exists_{W \in \Wfamily_{ucc}} \ W \subset V^{\dagger}
\label{exconditionexplicitely}
\end{gather}

\noindent
Intuitively, the condition (\ref{exconditionexplicitely}) means that the topology $\tau_{ucc}$ coincides with $\tau_{uc}$ on $\Ffamily$ (the inclusion $\tau_{ucc} \subset \tau_{uc}$ always holds so (\ref{exconditionexplicitely}) guarantees the reverse inclusion). It can be written more explicitely as 

\begin{gather}
\forall_{V \in \Vfamily} \ \exists_{\substack{D \Subset X \\ U\in\Vfamily}} \ \forall_{f,g \in \Ffamily} \ \bigg( \forall_{x \in D} \ g(x)\in U[f(x)] \ \Longrightarrow \ \forall_{x\in X} \ g(x) \in V[f(x)] \bigg)
\label{doesitlooklikebccondition}
\end{gather}

\noindent
If we suppose that $Y$ is a metric space, then $V[f(x)]$ becomes the ball $B(f(x),\eps)$ and $U[f(x)]$ becomes the ball $B(f(x),\delta)$. Hence (\ref{doesitlooklikebccondition}) turns into

\begin{gather}
\forall_{\eps > 0} \ \exists_{\substack{D - \text{compact} \\ \delta > 0}} \ \forall_{f,g \in \Ffamily} \ \bigg( \forall_{x \in D} \ g(x) \in B(f(x),\delta) \ \Longrightarrow \ \forall_{x\in X} \ g(x) \in B(f(x),\eps) \bigg)
\label{BCextcond}
\end{gather}

\noindent
which is the $BC(X,Y)$-extension property, studied in \cite{KrukowskiPrzeradzki}. Equipped with these tools, we are ready to prove the following theorem.

\begin{thm}(Arzel\`a-Ascoli for uniform convergence)\\
Let $(X,\tau_X)$ be a locally compact topological space, $(Y,\Vfamily)$ be a Hausdorff uniform space and let $C(X,Y)$ have the topology of uniform convergence. A subfamily $\Ffamily \subset C(X,Y)$ is relatively $\tau_{uc}$-compact if and only if 

\begin{description}
	\item[\hspace{0.4cm} (AA1)] $\Ffamily$ is pointwise relatively compact and equicontinuous
	\item[\hspace{0.4cm} (AA2)] $\Ffamily$ satisfies the extension property 
\end{description}  
\label{ArzelaAscoli}
\end{thm}
\begin{pro}
The '\textit{if}' part is easy. By \textbf{(AA1)} and theorem \ref{arzela-ascoliforuniformconvergenceoncompacta} we obtain that $\Ffamily$ is relatively $\tau_{ucc}$-compact. Consequently, \textbf{(AA2)} implies that the topology $\tau_{ucc}$ coincides with $\tau_{uc}$ on $\Ffamily$ and we are done.\\ 
We focus on '\textit{only if}' part. We already noted that $\tau_{uc}$ is stronger than $\tau_{ucc}$, which leads to $\Ffamily$ being relatively $\tau_{ucc}$-compact. Hence by theorem \ref{arzela-ascoliforuniformconvergenceoncompacta}, we obtain \textbf{(AA1)}. \\
Suppose that \textbf{(AA2)} is not satisfied, which means that there exists $V \in \Vfamily$ such that $W\backslash V^{\dagger} \neq \emptyset$ for every $W \in \Wfamily_{ucc}$. In particular, this means that 

\begin{equation}
\bigg\{(f,g) \in \Ffamily\times\Ffamily \ : \ \forall_{x\in D} \ (f(x),g(x)) \in V^{\frac{1}{9}}, \ \exists_{x_{\ast} \in X} \ (f(x_{\ast}),g(x_{\ast})) \not\in V\bigg\} \neq \emptyset
\label{negofextcond}
\end{equation}

\noindent
for every compact set $D$, where $V^{\frac{1}{9}}$ is a symmetric set such that $\underbrace{V^{\frac{1}{9}} \circ \ldots \circ V^{\frac{1}{9}}}_{9 \ \text{times}} \subset V$ and $A \circ B = \{(a,b) \ : \ \exists_{c} \ (a,c) \in B, \ (c,v) \in A\}$ as in \cite{Kelley} on page 176. The family 

\begin{gather}
\bigg\{g \in \Ffamily \ : \ \forall_{x \in X} \ (f(x),g(x)) \in V^{\frac{1}{9}} \bigg\}_{f \in \Ffamily}
\label{opencoverofF}
\end{gather}

\noindent
is a $\tau_{uc}$-open cover of $\overline{\Ffamily}$, the closure of $\Ffamily$ in $\tau_{uc}$. Indeed, if $\overline{f} \in \overline{\Ffamily}$, then 

$$(V^{\frac{1}{9}})^{\dagger}[\overline{f}] \cap \Ffamily \neq \emptyset \ \stackrel{(\ref{dagger})}{\Longleftrightarrow} \ \exists_{f \in \Ffamily} \ \forall_{x \in X} \ (f(x),\overline{f}(x)) \in V^{\frac{1}{9}}$$

\noindent
Since we assume that $\overline{\Ffamily}$ is $\tau_{uc}$-compact (and we aim to reach a contradiction), we can choose a finite subcover from (\ref{opencoverofF}), which means that there is a sequence $(f_k)_{k=1}^m \subset \Ffamily$ such that 

\begin{gather}
\forall_{g \in \overline{\Ffamily}} \ \exists_{k = 1,\ldots,m} \ \forall_{x \in X} \ (f_k(x),g(x)) \in V^{\frac{1}{9}}
\label{whatsequence}
\end{gather}

\noindent
Define (if it exists) an element $x_{kl} \in X$ for every $k,l = 1,\ldots,m$ such that 

$$(f_k(x_{kl}),f_l(x_{kl})) \not\in V^{\frac{1}{3}}$$

\noindent
and let $D$ be a compact set consisting of all $x_{kl}$'s. This set serves as a guard, watching whether each pair $f_k$ and $f_l$ 'drifts apart'. Its main task is the following implication

\begin{gather}
\forall_{x \in D} \ (f_k(x),f_l(x)) \in V^{\frac{1}{3}} \ \Longrightarrow \ \forall_{x \in X} \ (f_k(x),f_l(x)) \in V^{\frac{1}{3}}
\label{taskofD}
\end{gather}

\noindent
for every $k,l = 1,\ldots,m$. By (\ref{negofextcond}), we pick $(f_{\ast},g_{\ast})$ such that 

\begin{gather}
\forall_{x\in D} \ (f_{\ast}(x),g_{\ast}(x)) \in V^{\frac{1}{9}} \hspace{0.4cm} \text{and} \hspace{0.4cm} \exists_{x_{\ast} \in X} \ (f_{\ast}(x_{\ast}),g_{\ast}(x_{\ast})) \not\in V
\label{choiceoffgast}
\end{gather}

\noindent
Let $k_f, \ k_g \in \naturals$ be constants chosen as in (\ref{whatsequence}) for $f_{\ast}$ and $g_{\ast}$ respectively. We have 
	
$$\forall_{x \in D} \ (f_{k_f}(x),f_{\ast}(x)), \ (f_{\ast}(x),g_{\ast}(x)), \ (g_{\ast}(x),f_{k_g}(x)) \in V^{\frac{1}{9}} \ \Longrightarrow \ \forall_{x \in D} \ (f_{k_f}(x),f_{k_g}(x)) \in V^{\frac{1}{3}}$$

\noindent
By (\ref{taskofD}) we know that $(f_{k_f}(x),f_{k_g}(x)) \in V^{\frac{1}{3}}$ for $x \in X$. Finally, 

$$\forall_{x \in X} \ (f_{\ast}(x),f_{k_f}(x)), \ (f_{k_f}(x),f_{k_g}(x)), \ (f_{k_g}(x),g_{\ast}(x)) \in V^{\frac{1}{3}} \ \Longrightarrow \ \forall_{x \in X} \ (f_{\ast}(x),g_{\ast}(x)) \in V$$

\noindent
which is a contradiction with (\ref{choiceoffgast}). Hence, we conlcude that the extension property must hold. 
\end{pro}

\section{Applications}
\label{sectionapplications}

\noindent
In this final section of the paper, we present two possible applications of theorem \ref{ArzelaAscoli}. They were motivated by Schwartz's distribution theory, described in \cite{Schwartz} in chapter 2.

\subsection{Regular distributions}

\noindent
The space of distributions, denoted $\DST'(\reals^n)$, comes with the weak$^{\ast}$ topology (\cite{Schwartz} on page 94 or \cite{Duistermaat} on page 51). As such, it can be defined by psuedometrics 

$$\forall_{T,S \in \DST'(\reals^n)} \ p_{\vphi}(T,S) = |T(\vphi) - S(\vphi)|$$

\noindent
where $\vphi$ is a test function, i.e. compactly supported and infinitely differentiable function on $\reals^n$. The space of such functions is denoted by $\DST(\reals^n)$. It is easy to see that the family 

$$p_{\vphi}^{-1}([0,\eps)) = \bigg\{ (T,S) \in \DST'(\reals^n) \times \DST'(\reals^n) \ : \ |T(\vphi) - S(\vphi)| < \eps \bigg\}$$

\noindent 
for $\vphi \in \DST(\reals^n)$ and $\eps > 0$ is a subbase for the uniformity on $\DST'(\reals^n)$. Let $\Ffamily$ be a family of functions from $\reals \times \reals^n$ to $\reals$ such that

\begin{description}
	\item[\hspace{0.4cm} (M)] For every $t_{\ast} \in \reals$, the function $f(t_{\ast},\cdot)$ is measurable 
	\item[\hspace{0.4cm} (EQ)] For every $t_{\ast} \in \reals, \ \eps >0$ and compact set $D$, there exist $\delta > 0$ such that 
	
	 $$\forall_{\substack{t \in B(t_{\ast},\delta) \\ x\in D \\ f \in \Ffamily}} \ |f(t,x) - f(t_{\ast},x)| < \eps $$
	\item[\hspace{0.4cm} (BAB)] For every $t_{\ast} \in \reals$, the set $\{f(t_{\ast},x) \ : \ f \in \Ffamily, x \in \reals^n\}$ is bounded  
	\item[\hspace{0.4cm} (EP)] For every $\eps > 0$ and compact set $D \subset \reals^n$ there exists $R > 0$ such that for every $f,g \in \Ffamily$ we have 
	
	$$\sup_{t \in\reals\backslash B(0,R)} \ \int_D \ |f(t,x) - g(t,x)| \ dx < \eps$$
\end{description}

\noindent
The first two conditions are stronger than what is commonly known as \textit{Caratheodory conditions} (\cite{Aliprantis} on page 153). Condition \textbf{(EQ)} will turn out to be very useful, when proving equicontinuity. \textbf{(BAB)} will be used to establish pointwise relative compactness in the space $\DST'(\reals^n)$. To this end, we will employ Banach-Alaoglu-Bourbaki theorem, which can be found in \cite{Kothe} on page 248 or \cite{Meise} on page 264. Finally, with the use of \textbf{(EP)}, we will justify \textbf{(AA2)}. \\
Consider a family of distributions, defined by 

$$\forall_{\vphi \in \DST(\reals^n)} \ T_{f,t}(\vphi) := \int_{\reals^n} \ f(t,x)\vphi(x) \ dx$$

\noindent
where $f \in \Ffamily$ and moreover, denote $F_f : t \mapsto T_{f,t}$ (a distribution-valued function). We verify the assumptions of theorem \ref{ArzelaAscoli} in order to prove relative compactness of $F_{\Ffamily} := (F_f)_{f \in \Ffamily}$. \\
First of all, we check that $F_{\Ffamily}$ is equicontinuous, in particular $F_f \in C(\reals,\DST'(\reals^n))$ for every $f \in \Ffamily$. Condition (\ref{equi}) at fixed point $t_{\ast} \in \reals$ reads as 

\begin{gather}
\forall_{\substack{\vphi \in \DST(\reals^n) \\ \eps > 0}} \ \exists_{\delta > 0} \ \forall_{\substack{t \in B(t_{\ast},\delta) \\ f \in \Ffamily}} \ \bigg|\int_{\reals^n} \ (f(t,x) - f(t_{\ast},x))\vphi(x) \ dx\bigg| < \eps
\label{equidst}
\end{gather}

\noindent
For a fixed $\vphi \in \DST(\reals^n)$ and $\eps > 0$, we use \textbf{(EQ)} to choose $\delta > 0$ such that 

$$\forall_{\substack{t \in B(t_{\ast},\delta) \\ x\in \text{supp} \vphi \\ f \in \Ffamily}} \ |f(t,x) - f(t_{\ast},x)| < \frac{\eps}{\lambda(\text{supp} \ \vphi) \ \|\vphi\|}$$

\noindent
where $\lambda$ is the Lebesgue measure on $\reals^n$ and $\|\vphi\| :=  \sup_{x \in \reals^n} \ |\vphi(x)|$. Hence 

\begin{gather*}
\forall_{\substack{t \in B(t_{\ast},\delta) \\ f \in \Ffamily}} \ \bigg|\int_{\reals^n} \ (f(t,x) - f(t_{\ast},x))\vphi(x) \ dx\bigg| \leq \int_{\text{supp} \ \vphi} \ |f(t,x) - f(t_{\ast},x)||\vphi(x)| \ dx \\
< \eps \int_{\text{supp} \ \vphi} \ \frac{|\vphi(x)|}{\lambda(\text{supp} \ \vphi) \|\vphi\|} \ dx \leq \eps
\end{gather*}

\noindent
and we can conclude that $F_{\Ffamily}$ is equicontinuous. \\
In order to prove relative compactness of $F_{\Ffamily}$ at $t_{\ast}$, we use Banach-Alaoglu-Bourbaki theorem (\cite{Meise} on page 264 or \cite{Kothe} on page 248), which implies that for $U_0 \in \tau_{\DST(\reals^n)}$ and $r > 0$, the set  

\begin{gather}
\bigg\{T \in \DST'(\reals^n) \ : \ \sup_{\vphi \in U_0} \ |T(\vphi)| \leq r\bigg\}
\label{relcompdst}
\end{gather}

\noindent
is $\tau_{\DST'(\reals^n)}$-compact. Put 

$$U_0 = \bigg\{\vphi \in \DST(\reals^n) \ : \ \|\vphi\| < 1, \ \text{supp} \ \vphi \subset \overline{B}(0,1)\bigg\}$$

\noindent
which belongs to $\tau_{\DST(\reals^n)}$. Observe that 

$$\sup_{\vphi \in U_0} \ \bigg| \int_{\reals^n} \ f(t_{\ast},x)\vphi(x) \ dx\bigg| \leq \int_{\overline{B}(0,1)} \ |f(t_{\ast},x)| \sup_{\vphi \in U_0} \ \|\vphi\| \ dx \leq \lambda(\overline{B}(0,1)) \sup_{\substack{f \in \Ffamily\\ x\in\reals^n}} \ |f(t_{\ast},x)| $$  

\noindent
which proves that $F_{\Ffamily}(t_{\ast})$ is a subset of (\ref{relcompdst}) with $r = \lambda(\overline{B}(0,1)) \sup_{\substack{f \in \Ffamily\\ x\in\reals^n}} \ |f(t_{\ast},x)|$. Thus, we established pointwise relative compactness of $F_{\Ffamily}$. \\
Finally, we need to verify \textbf{(AA2)}. The uniformity of uniform convergence on $F_{\Ffamily} \subset C(\reals,\DST'(\reals^n))$ has a base of the form

\begin{gather*}
(p^{-1}_{\vphi}([0,\eps)))^{\dagger} = \bigg\{ (F_f,F_g) \in F_{\Ffamily} \times F_{\Ffamily} \ : \ \forall_{t \in \reals} \ (F_f(t),F_g(t)) \in p^{-1}_{\vphi}([0,\eps)) \bigg\} \\
= \bigg\{ (F_f,F_g) \in F_{\Ffamily} \times F_{\Ffamily} \ : \ \sup_{t \in \reals} \ |F_f(t)(\vphi) - F_g(t)(\vphi)| < \eps \bigg\} \\
= \bigg\{ (F_f,F_g) \in F_{\Ffamily} \times F_{\Ffamily} \ : \ \sup_{t \in \reals} \ \bigg|\int_{\reals^n} \ (f(t,x) - g(t,x)) \ \vphi(x) \ dx\bigg| < \eps \bigg\}
\end{gather*}

\noindent
Analogously, the uniformity of uniform convergence on compacta has a subbase of the form 

\begin{gather*}
W_{\vphi,R,\delta} = \bigg\{ (F_f,F_g) \in F_{\Ffamily} \times F_{\Ffamily} \ : \ \forall_{|t| \leq R} \ (F_f(t),F_g(t)) \in p^{-1}_{\vphi}([0,\delta)) \bigg\} \\
= \bigg\{ (F_f,F_g) \in F_{\Ffamily} \times F_{\Ffamily} \ : \ \sup_{t \in B(0,R)} \ \bigg|\int_{\reals^n} \ (f(t,x) - g(t,x)) \ \vphi(x) \ dx\bigg| < \delta \bigg\}
\end{gather*}
 
\noindent
Fix $\eps > 0$ and $\vphi \in \DST(\reals^n)$ and by \textbf{(EP)} choose $R > 0$ such that 

\begin{gather}
\forall_{f,g \in \Ffamily} \ \sup_{t \in\reals\backslash B(0,R)} \ \int_{\text{supp} \ \vphi} \ |f(t,x) - g(t,x)| \ dx < \frac{\eps}{2\|\vphi\|}
\label{epcondition}
\end{gather}

\noindent
Since

\begin{gather*}
\sup_{t \in \reals} \ \bigg| \int_{\reals^n} \ (f(t,x) - g(t,x)) \ \vphi(x) \ dx \bigg| \\
\leq \sup_{t \in B(0,R)} \ \bigg| \int_{\reals^n} \ (f(t,x) - g(t,x)) \ \vphi(x) \ dx \bigg| + \sup_{t \in \reals\backslash B(0,R)} \ \int_{\text{supp} \ \vphi} \ |f(t,x) - g(t,x)| \ dx \ \|\vphi\| \\
\stackrel{(\ref{epcondition})}{\leq} \sup_{t \in B(0,R)} \ \bigg| \int_{\reals^n} \ (f(t,x) - g(t,x)) \ \vphi(x) \ dx \bigg| + \frac{\eps}{2}
\end{gather*}

\noindent
we can conclude that for $\delta = \frac{\eps}{2}$ we have $W_{\vphi,R,\delta} \subset (p_{\vphi}^{-1}([0,\eps))^{\dagger}$, i.e. condition \textbf{(AA2)}. Thus we proved, using theorem \ref{ArzelaAscoli}, that $F_{\Ffamily} \subset C(\reals,\DST'(\reals^n))$ is relatively $\tau_{uc}$-compact.\\
In order to provide an example of a family $\Ffamily$ satisfying \textbf{(M)}, \textbf{(EQ)}, \textbf{(BAB)} and \textbf{(EP)}, consider 

$$\Ffamily = \bigg\{ a_ne^{-|g(t) - x|} \ : \ n\in\naturals\bigg\}$$  

\noindent
where $(a_n) \subset \reals$ is a bounded sequence and $g \in C(\reals,\reals)$ satisfies $\lim_{t\rightarrow\pm\infty} \ g(t) = \pm\infty$. Conditions \textbf{(M)} and \textbf{(BAB)} are obviously satisfied. For the remaining two conditions, fix $t_{\ast}\in\reals, \ \eps>0$, a compact set $D$ and set $a := \sup_{n\in\naturals} \ |a_n|$. By continuity, choose $\delta >0$ such that 

\begin{gather}
\forall_{t \in B(t_{\ast},\delta)} \ a\bigg(e^{|g(t)-g(t_{\ast})|} - 1\bigg) < \eps
\label{continuityofg}
\end{gather}

\noindent
Using the well-known estimate 

$$\forall_{u,v \geq 0} \ \bigg|e^{-u} - e^{-v}\bigg| \leq e^{-v}\bigg(e^{|u-v|} - 1\bigg)$$

\noindent
we can verify \textbf{(EQ)} as follows:

\begin{gather*}
\forall_{\substack{t \in B(t_{\ast},\delta) \\ x\in D\\ n\in\naturals}} \ \bigg|a_ne^{-|g(t)-x|} - a_ne^{-|g(t_{\ast})-x|}\bigg|
\leq ae^{-|g(t_{\ast})-x|}\bigg( e^{\big||g(t)-x| - |g(t_{\ast}) - x|\big|} - 1\bigg) \\
\leq a\bigg( e^{|g(t)-g(t_{\ast})|} - 1\bigg) \stackrel{(\ref{continuityofg})}{<} \eps
\end{gather*}

\noindent
It remains to prove \textbf{(EP)}. Since $\lim_{t\rightarrow\pm\infty} \ g(t) = \pm\infty$ we can choose $R>0$ such that 

\begin{gather}
\forall_{\substack{t \in \reals\backslash B(0,R) \\ x \in D}} \ |g(t) - x| > - \ln\bigg(\frac{\eps}{2a\lambda(D)}\bigg)
\label{limg}
\end{gather}

\noindent
For every $n,m \in \naturals$ we have 

$$\sup_{t\in B(0,R)} \ \int_D \ \bigg|a_ne^{-|g(t) -x|} - a_me^{-|g(t) - x|}\bigg| \ dx \leq 2a \sup_{t\in B(0,R)} \ \int_D \ e^{-|g(t) - x|} \ dx \stackrel{(\ref{limg})}{<} \eps$$

\noindent
which proves that $\Ffamily$ satisfies \textbf{(EP)}.

\section{Nonregular distributions - Dirac family}

\noindent
The second application of theorem \ref{ArzelaAscoli} concerns the following family of distributions :

$$\forall_{\vphi \in \DST(\reals^n)} \ T_{f,t}(\vphi) := \delta_{f(t)}$$

\noindent
where $\delta_p \in \DST'(\reals^n)$ is the Dirac delta function (\cite{Strichartz} on page 14) at point $p \in \reals^n$. We prove that if $\Ffamily \subset BC(\reals,\reals)$ is relatively compact, then the \textit{Dirac family} $\delta_{\Ffamily} := (\delta_f)_{f \in \Ffamily} \subset C(\reals,\DST'(\reals^n))$ is relatively $\tau_{uc}$-compact. \\ 
With the notation as before, we check the equicontinuity of $F_{\Ffamily}$. The condition (\ref{equi}) at a fixed point $t_{\ast} \in \reals$ reads as 

\begin{gather}
\forall_{\substack{\vphi \in \DST(\reals^n) \\ \eps > 0}} \ \exists_{r > 0} \ \forall_{\substack{t \in B(t_{\ast},r) \\ f \in \Ffamily}} \ |\vphi \circ f(t) - \vphi \circ f(t_{\ast})|
\label{equiFFfamily}
\end{gather}

\noindent
Fix $\vphi \in \DST(\reals^n)$ and $\eps > 0$. Since $\Ffamily$ is relatively compact, by uniform continuity of $\vphi$ on $\overline{\Ffamily}$ (closure of $\Ffamily$ in metric $d_{BC(\reals,\reals)}$) we obtain

\begin{gather}
\exists_{\rho > 0} \ \forall_{u,v \in \overline{\Ffamily}} \ |u - v|<\rho \ \Longrightarrow \ |\vphi(u) - \vphi(v)|<\eps
\label{unicontnonregulardist}
\end{gather}

\noindent
By equicontinuity of $\Ffamily$, we find that

\begin{gather}
\exists_{r > 0} \ \forall_{\substack{t \in B(t_{\ast},r) \\ f \in \Ffamily}} \ |f(t) - f(t_{\ast})| < \rho
\label{equinonregulardist}
\end{gather}

\noindent
By (\ref{unicontnonregulardist}) and (\ref{equinonregulardist}) we easily obtain (\ref{equiFFfamily}). \\
In order to prove pointwise relative compactness of $\delta_{\Ffamily}$ we need to show that 

\begin{gather}
\bigg\{\delta_{f(t_{\ast})} \ : \ f\in\Ffamily\bigg\}
\label{deltafamilyattast}
\end{gather} 

\noindent
is compact in $\DST'(\reals^n)$. Putting 

$$U_0 = \{\vphi \in \DST(\reals^n) \ : \ \|\vphi\| < 1\}$$

\noindent
we immediately see that $\sup_{\vphi \in U_0} \ |\vphi \circ f(t_{\ast})| \leq 1$ and consequently, (\ref{deltafamilyattast}) is $\tau_{\DST'(\reals^n)}$-compact.\\
Finally, as before, the family 

\begin{gather*}
(p^{-1}_{\vphi}([0,\eps)))^{\dagger} = \bigg\{ (\delta_f,\delta_g) \in \delta_{\Ffamily} \times \delta_{\Ffamily} \ : \ \sup_{t \in \reals} \ |\vphi \circ f(t) - \vphi \circ g(t)| < \eps \bigg\}
\end{gather*}

\noindent
is a base for uniformity of uniform convergence on $\DST'(\reals^n)$ and the uniformity of uniform convergence on compacta has a subbase of the form 

\begin{gather*}
W_{\vphi,R,r} = \bigg\{ (\delta_f,\delta_g) \in \delta_{\Ffamily} \times \delta_{\Ffamily} \ : \ \sup_{t \in B(0,R)} \ |\vphi \circ f(t) - \vphi \circ g(t)| < r \bigg\}
\end{gather*}
 
\noindent
It becomes evident, that since $\Ffamily$ satisfies $BC$-extension condition (\ref{BCextcond}), then for fixed $\vphi \in \DST(\reals^n)$ and $\eps> 0$, there exist $R,r > 0$ such that $W_{\vphi,R,r} \subset (p_{\vphi}^{-1}([0,\eps)))^{\dagger}$, meaning \textbf{(AA2)} is satisfied. We conclude, via theorem \ref{ArzelaAscoli}, that $\delta_{\Ffamily}$ is relatively $\tau_{uc}$-compact.

\Addresses
\end{document}